%% file: Paper2.tex
\documentclass[a4paper,12pt]{article}

\usepackage[left=2.5cm, right=2.5cm, top=2.5cm, bottom=2.5cm]{geometry}

\usepackage{amsmath}
\usepackage{amssymb}
\usepackage{amsthm} 
\usepackage{mathtools}
\usepackage[T1]{fontenc}
\usepackage[utf8]{inputenc}
\usepackage{pdfpages}

\usepackage[sorting=nyt, style = ieee]{biblatex}

\usepackage{hyperref}
\usepackage{cleveref}
\usepackage{caption}
\usepackage{subcaption}

\newtheorem{theorem}{Theorem}[section]
\newtheorem{corollary}[theorem]{Corollary}

\def\norm#1{\lVert #1 \rVert}

\bibliography{eigen.bib}

\newenvironment{@abs}[1]{%
  \vspace{4pt}%\footnotesize
  \parindent 15pt {\bfseries #1. }\ignorespaces
     }
     {\par\vspace{7pt}}
     
\renewenvironment{abstract}{\begin{@abs}{{Abstract}}}{\end{@abs}}
\newenvironment{keywords}{\begin{@abs}{Key words}}{\end{@abs}}
\newenvironment{AMS}{\begin{@abs}{AMS subject classifications}}{\end{@abs}}

\title{Some notes concerning preconditioning
  of linear parabolic optimal control problems}
\author{ Luise Blank \footnotemark[2]
}
\date{}

\begin{document}

\maketitle

\renewcommand{\thefootnote}{\fnsymbol{footnote}}

\footnotetext[2]{Department of Mathematics,University of Regensburg, 93040 Regensburg, Germany.}

\input{eigen2ArXiv}

\printbibliography

\end{document}

%% file: eigen2ArXiv.tex
\begin{abstract}
In this paper we study the conditioning of optimal control problems constrained by  linear parabolic equations with Neumann boundary conditions.
While we concentrate on a given end-time target function the results hold also when the target function is given over the whole time horizon.
When implicit time discretization and conforming finite elements in space are employed we show that the reduced problem formulation has condition numbers which are bounded independently of the discretization level in arbitrary space dimension.
In addition we propose for the all-at-once approach, i.e. for the first-order conditions of the unreduced system a preconditioner based on work by Greif and Schötzau, which provides also bounds on the eigenvalue distribution independently of the discretization level.
Numerical experiments demonstrate the obtained results and the efficiency of the suggested preconditioners.
\end{abstract}



%%%%%%%%%%%%%%%%%%%%%%%% 
\section{Introduction}
The development of fast solvers for optimizing nonlinear parabolic problems hinges in many cases on fast solvers for the arising linear systems. It is not clear whether the elimination of the state given by the uniquely solvable differential equation for a given control is more efficient than considering the full system,
even for minimizing linear parabolic problems.
One point is the conditioning of the arising linear
systems. 
Considering the first-order conditions for the full system leads to a saddle point system  involving
the control, the state and the adjoint and
therefore is also
called the all-at-once approach. For nonlinear constraints one has to keep in mind that a globally convergent solver has to be employed.
If an optimization solver for nonlinear equality constrained problems is applied, as e.g. an SQP-method, also changing linear saddle point problems have to be solved in each iteration.
For the reduced optimization problem solvers for unconstrained problems can be applied.
Then, in general the discretized PDEs have to be solved with a rather high accuracy in each iteration.
There is vast literature concerning preconditioning saddle point problems in general (see
\cite{BenziGolubLiesen2005}, \cite{Wathen2015}  
 and references therein). 
For saddle point problems arising from optimization problems with elliptic PDEs also many papers exist
(see \cite{ReesDollarWathen2010,SchoeberlZulehner2007,SchoeberlSimonZulehner2011} and references therein).
A survey of preconditioners for Krylov subspace methods categorized by application area is given in  \cite{PearsonPestana2020}.
It includes 
also a summary of more recent results concerning preconditioners for optimization problems constrained by PDEs. 
In 
\cite{MathewSarkisSchaerer2010, PearsonStollWathen2012}
control problems constrained by the heat equation and with a target function given over the whole time horizon (tracking function) are studied.
In \cite{PearsonStoll2013} preconditioners for problems with more general reaction diffusion equations with or without control constraints are presented and numerically analyzed.
In all three papers a fixed time step and a spatial mesh constant in time are employed.
They develop preconditioners in particular for the all-in-once approach.
In \cite{SchielaUlbrich2014}
operator preconditioning is studied where the focus lies on control or state constraints. The underlying linear  PDEs can be time-dependent involving a coercive elliptic operator. Due to the inequality constraints the system is reduced to the adjoint and the state or to the adjoint only.
Numerical experiments  are given for tracking type problems for the Laplace equations with inequality constraints.

\medskip
In our paper 
we concentrate on control problems with
linear parabolic state equations 
with Neumann boundary conditions
where the elliptic part is not necessarily coercive.
We focus on cost functionals with a given end-time target function.
However, the results for the reduced problem hold also for a given tracking function. 
The time steps and the spatial mesh may vary over time.
The paper is structured as follows. 
In the next section we introduce the problem formulation, its discretization, which uses dG(0) in time, i.e. an
implicit time discretization, and conforming finite elements in space, and  set the notations.
For the reduced optimality system
we derive mesh independent bounds for the condition numbers in Section 3. Here, a-priori estimates known in the  analysis of PDEs are applied. Due to the use of conforming Galerkin discretization, the estimates hold also for the fully discretized problem.
Numerical evidence is shown as well.
In the following section we present 
for the all-at-once approach a preconditioner
providing mesh  independent bounds on the eigenvalue distribution.
Here, the resulting saddle point problem has a highly singular 1-1-block reflecting that the target function is only given at the final time. The employed preconditioner is based on the work of Greif and Schötzau
\cite{GreifSchoetzau2006}.
Numerical experiments confirm the analytical results.

%%%%%%%%%%%%%%%%%%%%%%%%
%%%%%%%%%%%%%%%%%%%%%%%%
\section{Problem formulation and time discretization}
The considered optimal control problems of a linear parabolic partial differential equation with Neumann boundary conditions
read as
\begin{eqnarray}
  \label{cost1}
 && \min_{u\in L^2(Q)}  J(u,y) \coloneqq \dfrac{1}{2} \norm{y(T) - y_\Omega}^2_{L^2(\Omega)} + \dfrac{\lambda}{2} \norm{u}^2_{L^2(Q)}\\
  \mbox{or}\nonumber &&\\
  \label{cost2}
 && \min_{u\in L^2(Q)}  J_Q(u,y) \coloneqq \dfrac{1}{2} \norm{y- y_Q}^2_{L^2(Q)} + \dfrac{\lambda}{2} \norm{u}^2_{L^2(Q)}
  \end{eqnarray}
subject to the state equation
\begin{equation}
\label{pde}
\begin{aligned}
  (\partial_t y , \eta)_{L^2(Q)}+(\alpha \nabla y, \nabla \eta)_{L^2(Q)}
+(c y, \eta)_{L^2(Q)}
  &=  (u,\eta) _{L^2(Q)}
     \;  \forall \eta \in L^2(0, T; H^1(\Omega))\\
     y(0)&=y_0     \quad \mbox{ in } \Omega
     \end{aligned}
\end{equation}
where $\Omega $ denotes the space domain in $\mathbb{R}^d$, $Q:=  [0,T]\times \Omega$ is the space-time cylinder and where 
the coefficient functions fulfill $\alpha\geq 0 $ and $\alpha, c\in L^\infty (Q) $.
The state equation is formulated in the weak sense and includes therefore the homogeneous Neumann boundary condition.
While with the cost functional $J$ one tries to find a distributed control $u$ such that the state $y$ is close to a given target function $y_\Omega \in L^2(\Omega)$ at the end-time $T$,
the functional $J_Q$  aims to obtain a minimal distance of $y$ to a target function $y_Q \in L^2(Q)$ over the whole time horizon. Since the former leads to more complicated structures we focus on $J$.
The results of Section 3  carry over to $J_Q$. They also hold for the state equation with $H^1_0(\Omega )$ instead of $H^1(\Omega )$, i.e. with Dirichlet boundary conditions. Moreover, the minimization may be only over $u \in L^2(\bar Q) $ with $\bar Q \subseteq Q$ with nonzero measure, or the cost functional involves only $ \norm{y(T) - y_{\bar \Omega}}^2_{L^2(\bar \Omega)} $ with $\bar \Omega  \subseteq \Omega$ with nonzero measure.

Let us state briefly the optimality conditions for $J$
(see e.g. \cite{Lions71, HinzePinnauUlbrichUlbrich2009, Troeltzsch2010}). They are given by state equation \eqref{pde}, by the adjoint equation 
\begin{equation}
\label{adjoint}
\begin{aligned}
-  (\eta,\partial_t p)_{L^2(Q)} +  (\alpha \nabla \eta, \nabla p)_{L^2(Q)} + ( c \eta,p)_{L^2(Q)} &= 0 \qquad \forall \eta \in L^2(0, T;H^1(\Omega))\\
p(T) &= y(T)-y_\Omega   \quad\text{in } \Omega
.
\end{aligned}
\end{equation}
and by
\begin{equation}
  \label{KKT1}
  \lambda {u} + {p} =0 \quad \text{ a.e. in } Q. 
\end{equation}
For the reduced cost functional $j(u):=J(y(u), u)$
the gradient is also determined by
$ \nabla{j}({u})  = \lambda {u} + {p}$.
These equations are now discretized in time in such a way that
the arising system coincide with the system when
the minimization problem is discretized first and
then the optimality system is set up.
Therefore we use the finite element discretization
dG(0) which is equivalent to the implicit Euler discretization \cite{NeittaanmakiTiba94}. This means, given a time mesh $0=t_0<t_1<\ldots <t_N=T$  with time steps $\tau_n:=t_n-t_{n-1}$
the involved time discretized functions
$u_\tau, y_\tau,p_\tau$ 
are considered to be constant in time on $I_n:=(t_{n-1},t_n)$ with values $u_n \in L^2(\Omega)$, respectively
$y_n, p_n \in H^1(\Omega)$ for $n=1,\ldots,N $.
Having the view point of an implicit Euler discretization we identify $u_n=u_\tau(t_n), y_n=y_\tau(t_n)$ while
$p_n=p_\tau(t_{n-1})$.
In the remainder of the  paper $(.,.)$ and $ ||\cdot ||$ denote the ${L^2(\Omega)}$-scalar-product and
${L^2(\Omega)}$-norm.
We obtain the time discretized problem
\begin{equation}
\label{discrete_problem1}
\min
 {J}({y_\tau},{u_\tau})
= \dfrac{1}{2}\|y_N-y_\Omega\|^2 + \dfrac{\lambda}{2}\sum_{n=1}^N \tau_n \|u_{n}\|^2
\end{equation}
\begin{equation}
\begin{aligned}
  \text{s.t.} \quad
  \label{scheme_state1}
\tfrac1{\tau_{n}}(y_{n}-y_{n-1}, \varphi_{}) +   (\alpha_n \nabla y_{n}, \nabla \varphi_{}) +
(c_n y_{n}, \varphi_{}) &=
(u_{n}, \varphi_{})
\quad \forall \varphi \in H^1(\Omega), n = 1,\ldots,N 
\end{aligned}
\end{equation}
where $y_0 \in H^1(\Omega)$ is the given initial value.

The discrete adjoint equation reads as:
\begin{equation}	
\begin{aligned}
\label{eq:disc_adjoint}
\tfrac1{\tau_{j}}(\varphi_{},p_{j}-p_{j+1}) +(\alpha_j \nabla \varphi, \nabla p_{j}) +  (c_j \varphi ,p_{j}) &=
0 
\quad \forall \varphi \in H^1(\Omega), j=N,\ldots, 1 
\end{aligned}
\end{equation}
with $p_{N+1}:=y_N-y_\Omega$.
The optimality condition and the gradient of
$j_\tau (u_\tau):=J(u_\tau, y_\tau) $ 
is given as in the continuous version by 
\begin{align}\label{gradjtau}
  \nabla{j_\tau }({u_\tau}) &=
                              \lambda { u_\tau} + {p_\tau} =0
                              \quad
                              \text{ a.e. in } \Omega .
\end{align}
The fully discretized problem using conforming finite elements  $S^{(j)}_h\subset H^1(\Omega )$  in space
\cite{Ciarlet2002, Bartels2015} 
for $u_j, y_j$ and $p_j$  on each Interval $I_j$ is given by the same formula as in the time discretized case employing test functions $\phi \in S^{(j)}_h$.
The spatial mesh may be nonuniform and may differ in each time step. We denote the corresponding functions by $u_{\tau, h}$ and so on.
Fixing the basis for the discretization spaces,
i.e.  let $\phi^{(j)}_l$ be the basis functions for the finite element space on the $j-$th time interval $I_j$,
we use the notation $\bar u$ etc. for the corresponding coefficients.

We define the following matrices:  The symmetric positive
definite mass matrix on $I_j$ is given by
 \begin{eqnarray}
   M_j:=\left( (\phi^{(j)}_l ,\phi^{(j)}_k) \right) _{k.l}
\;    , \quad \mbox{and} \quad
   {\cal M}:=blkdiag(M_j)_{j=1,\ldots,N}
   \; , \label{MatrixMj}
 \end{eqnarray}   
 
 the matrix  $  \cal K$ is  the discretized operator for $y$,
 i.e.
 \begin{eqnarray}
   \label{MatrixK}
 {\cal K} &:=& blkdiag(K_j )_{j=1,\ldots , N}
          - blkdiag(
 \tfrac1{\tau_j} \left( (\phi^{(j-1)}_l ,\phi^{(j)}_k) \right) _{k.l}
               ,-1)_{j=2,\ldots , N}\\
\label{MatrixKj}
\mbox{with }  K_j&:=&\tfrac1{\tau_j} M_j
 +\left( (\alpha_j \nabla \phi^{(j)}_l , \nabla \phi^{(j)}_k) 
 +  (c_j \phi^{(j)}_l ,\phi^{(j)}_k) 
   \right) _{k.l}
\end{eqnarray}
and let  $D:= blkdiag(\tau_j I)_{j=1,\ldots,N}.$
Then  the unreduced optimality system is given by
 \begin{eqnarray}
   {\cal K} \bar y -{\cal M}\bar u &=&(\tfrac1{\tau_1} M_{1,0} \bar y_0, 0,\ldots,0)^T\\
 \tfrac1{\tau_N} e_N e_N^T \otimes M_N \bar y
 - D^{-1} {\cal K}^{T} D \bar p 
 &=& (0,\ldots,0,\tfrac1{\tau_N} M_{N,\Omega} \bar y_\Omega )^T\\
 \lambda \bar u +\bar p &=&0.
 \end{eqnarray}
 Here $M_{1,0} \bar y_0$ and $ M_{N,\Omega}\bar y_\Omega$ are the vectors 
 given by $((y_0 ,\phi^{(1)}_k))_k$ and
 $((\phi^{(N)}_l ,y_\Omega))_l$ or their approximations.
 Written in a symmetric way it reads as 
\begin{equation}
	\label{KKTsym}
\underbrace{\left(   \begin{array}{ccc}
      e_N e_N^T \otimes M_N & 0 & {{\cal K }}^T D\\
0  & \lambda  {D\cal M} & - {\cal M}D
\\
{D\cal K} & - {D\cal M}  & 0 
\end{array}\right)
 }_{=: \cal T}               
 \left( \begin{array}{c}
\overline y \\ \overline u \\ - \overline p 
\end{array}\right)
    = 
 \left(     \begin{array}{c}
 \left(     \begin{array}{c} 0\\  M_{N,\Omega} \overline y_\Omega\end{array}\right)
      \\
      0
      \\
  \left(      \begin{array}{c}  
  	M_{1,0} \overline y_0\\ 0  \end{array}\right)
    \end{array}\right).
\end{equation}
First discretizing then optimizing leads to a Lagrange-multiplier $\zeta = D \bar p $ and the optimality system
\begin{equation}
	\label{KKTdisc}
	\underbrace{\left(   \begin{array}{ccc}
			e_N e_N^T \otimes M_N & 0 & {{\cal K }}^T \\
			0  & \lambda  {D\cal M} & - {\cal M}
			\\
			{\cal K} & - {\cal M}  & 0 
		\end{array}\right)
	}_{=: {\cal T}_{disc}}               
	\left( \begin{array}{c}
		\overline y \\ \overline u \\ - \zeta
	\end{array}\right)
	= 
	\left(     \begin{array}{c}
		\left(     \begin{array}{c} 0\\  M_{N,\Omega} \overline y_\Omega\end{array}\right)
		\\
		0
		\\
		\left(      \begin{array}{c}  
		\tfrac1{\tau_1}	M_{1,0} \overline y_0\\ 0  \end{array}\right)
	\end{array}\right).
\end{equation}

%%%%%%%%%%%%%%%%%%%%%%%%
%%%%%%%%%%%%%%%%%%%%%%%%
\section{Eigenvalue estimates for the reduced problem formulation}
The optimization problem may be solved numerically by considering the problem reduced to the control $u$. Then  the  affine system
\begin{equation}
  \lambda u_\tau+p_\tau=0
\end{equation}
is solved iteratively, where 
it is assumed that the discretized state and adjoint equations are solved for a given control $u_\tau$ exactly.
The state $y_\tau(u_\tau,y_0)$ depends on the right hand side $u_\tau$  and the given initial data $y_0$.
For the adjoint $p_\tau(y_N-y_\Omega)$ the right hand side is $0$ and the final value is given by $y_N$ and the given data $y_\Omega$.
Due to linearity we have
$p_\tau(y_N(u_\tau,y_0)-y_\Omega)
=p_\tau(y_N(u_\tau,0)) + p_\tau(y_N(0,y_0)-y_\Omega)
$.
Hence one solves the linear system
\begin{equation}
  F_\tau  u_\tau := \lambda u_\tau +p_\tau (y_N(u_\tau ,0))= - p_\tau(y_N(0,y_0)-y_\Omega)
  \label{Ftau}
 \end{equation}
where $p_\tau(y_N(u_\tau,0))$ solves the adjoint equation with $p_{N+1}=y_N(u_\tau,0)$ and $y_\tau(u_\tau,0)$ solves the state equation with initial data $y_0=0$.

\medskip
For the eigenvalue estimates we use a-priori estimates for $y_\tau$ and $p_\tau$, which we briefly show here for completeness although it can be found in the literature (see e.g. \cite{Evans2010}).

With $\tfrac{1}{2}(a^2-b^2) \leq (a-b)a$,
 testing the state equation \eqref{scheme_state1} with $ y_n$,
employing the scaled Young's inequality  with a scaling parameter $\epsilon >0$ and
using the existence
of a constant $c_0 \in \mathbb{R}$ with  $c_0 \leq c_n(x) $   a.e. in $\Omega $  for all $N$ and $n=1,\ldots N$
due to $c \in L^\infty(Q)$
gives
{\begin{align*}
  \tfrac{1}{2}\left(\norm{y_n}^2 - \norm{y_{n-1}}^2\right)
&\leq \left(y_n - y_{n-1},  y_n\right)
+ \tau_n \left(\alpha_n \nabla y_{n}, \nabla y_n \right) 
   \\&
   = \tau_n \left(u_n, y_n\right)
    - \tau_n  \left( c_n y_{n},  y_n\right)
   \leq \tfrac{\tau_n}{2\epsilon} \norm{ u_n}^2 + \tfrac{\tau_n}{2}(\epsilon-2c_0) \norm{ y_n}^2
\\
  \Rightarrow
  \norm{y_J}^2 \leq  \tfrac1{1-\tau_J(\epsilon-2c_0)}
    &  \left( \norm{y_0}^2+
      \tfrac{1}{\epsilon} \norm{u_\tau}^2_{L^2(Q)} 
      + (\epsilon-2c_0) \sum_{n=1}^{J-1} \tau_n \norm{ y_n}^2\right)
\end{align*}
}
if $ \tau_J(\epsilon-2c_0) <1$.
Choosing $\epsilon \geq 2c_0$ we can apply
the discrete Gronwall Lemma
\begin{eqnarray}
\label{eqgron1}
 \norm{y_J}^2 &\leq &\left(\norm{y_0}^2 + \tfrac{1}{\epsilon}  \norm{u_\tau}^2_{L^2(Q)}  \right)
 \tfrac1{1-\tau_J(\epsilon-2c_0)}
                      \exp\left( \frac{\epsilon-2c_0}{1-\tau_J(\epsilon-2c_0)} T\right) \\
  \label{eqgron2}
              &\leq &\gamma \left(\epsilon \norm{y_0}^2 + \norm{u_\tau}^2_{L^2(Q)}  \right)
\end{eqnarray}
where $\gamma $ depends only on $c_0$ and $T$
but  is independent of $\tau_J$ for small enough   $\tau_J$, i.e. when for some
chosen
$\epsilon \geq 2c_0$ with $\epsilon >0  $ and
$\delta <1$ it holds
$ \tau_J(\epsilon-2c_0)\leq \delta$.
To see the influence of $c_0$ and $T$ on the bound $\gamma$ let us give possible choices of $\epsilon$.
If $2c_0> -\tfrac1T$ one can choose e.g.
$\epsilon= \tfrac1T+2c_0$.
Then we have 
$\gamma \leq 3 \tfrac{T}{1+2c_0T} $
for all $\tau_J \leq 0.001 T$ .
Otherwise the choice $\epsilon =\tfrac1T$ leads to the bound
$\gamma \leq 3 T \exp(- 2.002 c_0T) $
for 
all  $\tau_J\leq 0.001 T/ (1-2c_0T) $ .

If the input $u_\tau \equiv 0$ then one obtains
correspondingly for $ - \tau_{J} 2 c_0<1$
\begin{eqnarray}
  \norm{y_J}^2 &\leq &\norm{y_0}^2
                       \left\{
 \begin{array}{l l}                     1
                       & \text{ if } c_0 \geq0 \\
 \tfrac1{1+\tau_J2c_0}
                       \exp\left( \frac{-2c_0}{1+\tau_J2c_0} T\right)
&\text{ if } c_0 < 0 
 \end{array}
   \right.
  .
  \label{gamma}
\end{eqnarray}
Rewriting the equation for the adjoint $p_\tau$,
which is determined backward in time,
as a forward problem one can apply the above inequality. Let for $\tau_{max}:=\max\{\tau_1,\ldots,\tau_N \}$ hold $- \tau_{max} 2 c_0<1$.   Then we obtain the estimate
\begin{eqnarray}
  \qquad
  ||p_\tau||_{L^2(Q)}^2
  \leq
T \norm{y_N}^2 \tfrac1{1+\tau_{max}2\min\{c_0,0\}}
 \exp\left( \frac{-2\min\{c_0,0\}}{1+\tau_{max}2c_0} T\right)
 \leq  T \gamma_2 \norm{y_N}^2
  \label{eqp}
\end{eqnarray}
where the last inequality holds with some $\gamma_2>0$ for all small enough $\tau_{max}$. 

\medskip
For the fully discretized system using conforming finite elements in space the same steps and estimates hold true
independently of the spatial mesh.
We denote by $F_{\tau,h} $ the corresponding
linear map given by
\begin{equation}
    F_{\tau,h}  u_{\tau,h} := \lambda u_{\tau,h} +p_{\tau,h} (y_{\tau,h}|_T(u_{\tau,h} ,0)) 
\end{equation}
and with $ \bar F_{\tau,h }$
its matrix representation
using the basis $(\phi_l^{(j)})$.
Hence we have
\begin{eqnarray}
   \bar F_{\tau,h }\bar u =
  \lambda \bar u +\bar p (\bar y_N (\bar u,0) ) 
 =        (
  \lambda I + D^{-1}{\cal K}^{-T}(e_N e_N^T \otimes M_N)
  {\cal K}^{-1} {\cal M}
  )\bar u  \; .
  \label{27}
\end{eqnarray} 
The estimates provide us with the following eigenvalue inclusion, where we use the notation
$\eta(A)$ for eigenvalues of a matrix $A$.
\begin{theorem}
 There exist $\gamma , \bar \tau > 0$, such that
 for  all eigenvalues
 $\eta(F_{\tau})$, $\eta(F_{\tau,h})$,  $\eta(\bar F_{\tau,h})$
 it holds in any space dimension
 \begin{equation}
   \label{eigenF}
  \eta(F_{\tau}),  \eta(\bar F_{\tau,h}) = \eta(F_{\tau,h})\in [\lambda, \lambda +\gamma]
  \qquad \forall \tau_N\leq \bar \tau
  \text{ and } \forall h> 0.
\end{equation}
Furthermore, at most
$rank (e_N e_N^T \otimes M_N )=\dim (S_h^{(N)}) $
eigenvalues of  $\bar F_{\tau,h}$ as well as of $F_{\tau,h}$
do not equal $\lambda$.
\end{theorem}
\begin{proof}
Without restriction we consider $F_\tau$. For a given eigenfunction $u_\tau$ it holds
  $ 
  \eta(F_\tau)  ||u_\tau||^2_{L^2(Q)}
  =(u_\tau, F_\tau u_\tau)_{L^2(Q)}
  = \lambda ||u_\tau||^2_{L^2(Q)} +
  (u_\tau, p_\tau(y_N(u_\tau,0))_{L^2(Q)}
  $.
  Furthermore the derivation of the adjoint
  and \eqref{eqgron2} using $y_0=0$ guarantees the existence of $\gamma \geq 0$ and $\bar \tau >0$ with
  $ (u_\tau, p_\tau)_{L^2(Q)}=
  (y_N,y_N)_{L^2(\Omega)}
  \leq \gamma ||u_\tau||^2_{L^2(Q)}$
  for all $\tau_N\leq \bar \tau$ .
  Hence it holds
  \begin{eqnarray*}
  \lambda ||u_\tau||^2_{L^2(Q)}
 & \leq& \lambda ||u_\tau||^2_{L^2(Q)} + (y_N,y_N)_{L^2(\Omega)}
         =\eta(F_\tau)  ||u_\tau||^2_{L^2(Q)}
\\  &\leq& \lambda ||u_\tau||^2_{L^2(Q)} + \gamma ||u_\tau||^2_{L^2(Q)}.
\end{eqnarray*}
\hfill \end{proof}

The constant in the eigenvalue estimate may be sharpened by using \eqref{eqp}
\begin{eqnarray}
  (u_\tau,p_\tau)_{L^2(Q)}\leq ||u_\tau||_{L^2(Q)}  ||p_\tau||_{L^2(Q)}\leq ||u_\tau||^2_{L^2(Q)}
  \sqrt{T \gamma \gamma_2}.
\end{eqnarray}
The matrix $\bar F_{\tau,h}$ is symmetric with respect to the $L^2(Q)$-scalar-product discretized with the given basis functions, i.e. with respect to $D\cal M$.
Hence this scalar product has to be used when applying e.g. the CG-method to
$\lambda \bar u +\bar p =rhs$.
Then due to \eqref{eigenF}
the number of CG-iterations are bounded independently of the discretization level.

Considering the symmetrically formulated optimality system \eqref{KKTdisc}
and reducing it to  $\bar u$ we obtain
\begin{equation}
 {\cal  R}:= \lambda D{\cal M} +  
 {\cal M}{\cal K} ^{-T} (e_N e_N^T \otimes M_N) {\cal K} ^{-1} {\cal M} \; .
 \label{R}
\end{equation}
With $({D\cal M})^{-1} {\cal R}= \bar F_{\tau,h}$ the following result immediately holds:
\begin{corollary}
  The matrix $\cal R $ is symmetric positive definite
  and there exist $\gamma , \bar \tau > 0$ such that for all $\tau$ with $0< \tau_N \leq \bar \tau$ and for all $h>0$
the eigenvalues of preconditioned matrix
$({D \cal M})^{-1} \cal R$ lie in
{$[\lambda,\lambda +\gamma] $.}
Furthermore, at most $\dim (S_h^{(N)})$ eigenvalues do not equal $\lambda $.
\end{corollary}

The above considerations carry over to the discretized optimization problem \eqref{cost2} where the target function is given over the whole time-horizon
and also to the control problem when Dirichlet data are given. 
\begin{corollary}
The matrix ${\cal R }_Q:= \lambda D{\cal M} +  
  {\cal M}{\cal K} ^{-T} (D{\cal M}) {\cal K} ^{-1} {\cal M}$ is symmetric positive definite
  and there exist $\tilde \gamma , \bar \tau > 0$ such that for all $\tau$ with $0< \max\{\tau_i \} \leq \bar \tau$ and for all $h>0$
the eigenvalues of preconditioned matrix
$({D \cal M})^{-1} {\cal R}_Q$ lie in
{$[\lambda,\lambda +\tilde \gamma] $.}
\end{corollary}

For the specific case of Dirichlet boundary conditions and without the lower order term
$ (c y, \eta)_{L^2(Q)}$ in the state equation
this result is also a consequence of Theorem 3.2 in \cite{MathewSarkisSchaerer2010} if
uniform meshes in space and time are used.

Given above estimates the PCG-method applied to $\cal R$ with preconditioner ${D \cal M}$ yields iteration numbers
bounded independently of the discretization level.
Let us mention that it is cheaper to apply the CG-method to 
$\bar F_{\tau,h}$ with the $D{\cal M}$-scalar product. 
Here, one can reduce additionally the cost of the method by one application of $D{\cal M}$ by introducing an auxiliary variable like in the PCG-method.
%%%%%%%%%%%%%%%%%%%%%%%%%%%%%%
\subsection{Numerical evidence}\label{evidencePcg}
The implementation is done in MATLAB using the toolbox for partial differential equations. As spatial domain we choose $\Omega =(0,1)^d$.
We employ equidistant time steps and the same spatial mesh in each time
step as well as constant functions $\alpha $ and $c$.
Hence the mass matrices $M_j$ and stiffness matrices
$K_j=\tfrac1{\tau_j} M_j
 +\left( (\alpha_j \nabla \phi^{(j)}_l , \nabla \phi^{(j)}_k) 
 +  (c_j \phi^{(j)}_l ,\phi^{(j)}_k) 
 \right) _{k.l}
 $
 given in \eqref{MatrixMj} and \eqref{MatrixKj}
stay constant in time, i.e with respect to $j$.
The Cholesky decompositions of these two matrices are performed, and then they are employed in each iteration in each time step for the state as well as for  the adjoint equation. 

\begin{figure}[htbp]
\begin{center}
  \includegraphics[width=.4\textwidth]{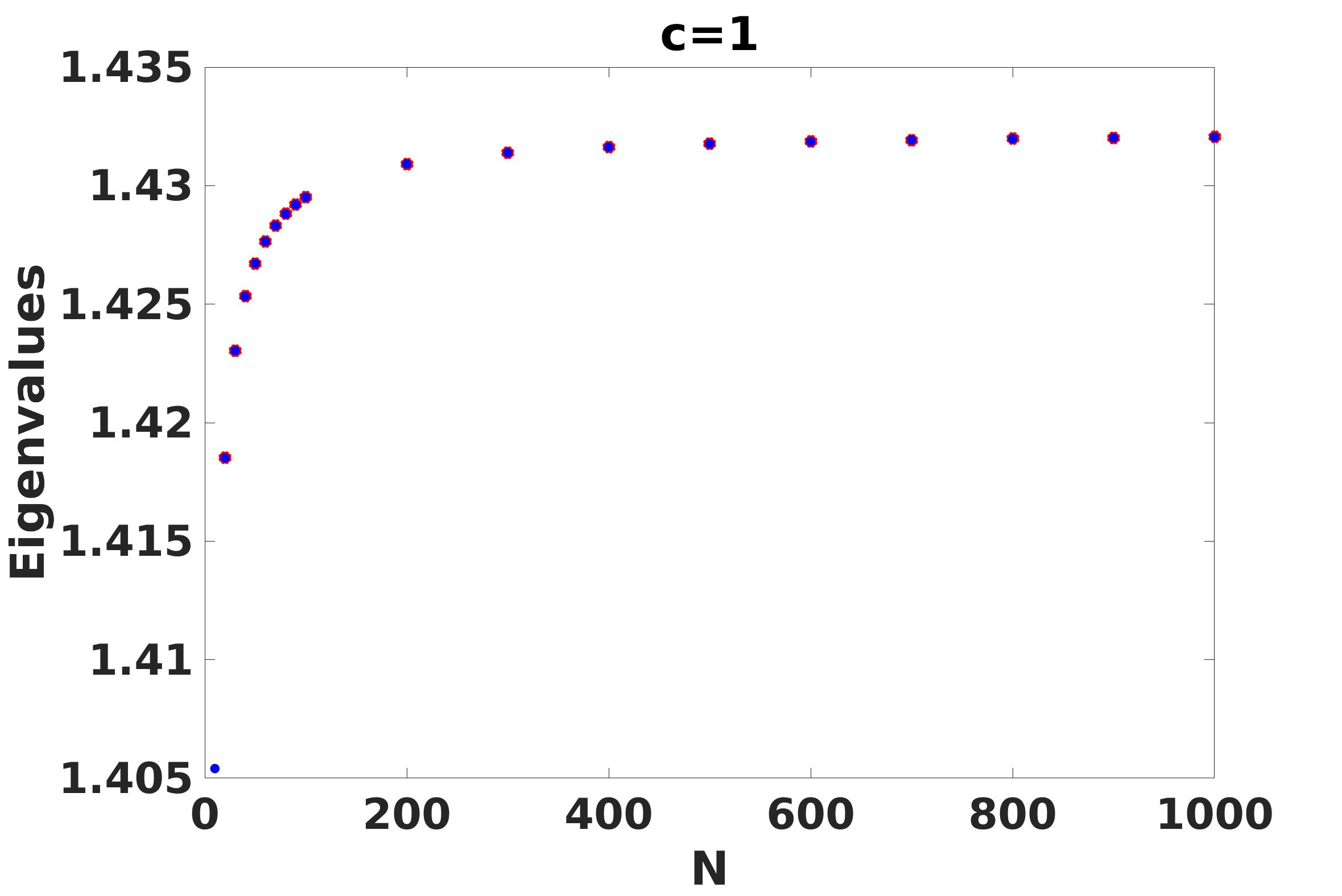}~~~~
\includegraphics[width=.4\textwidth]{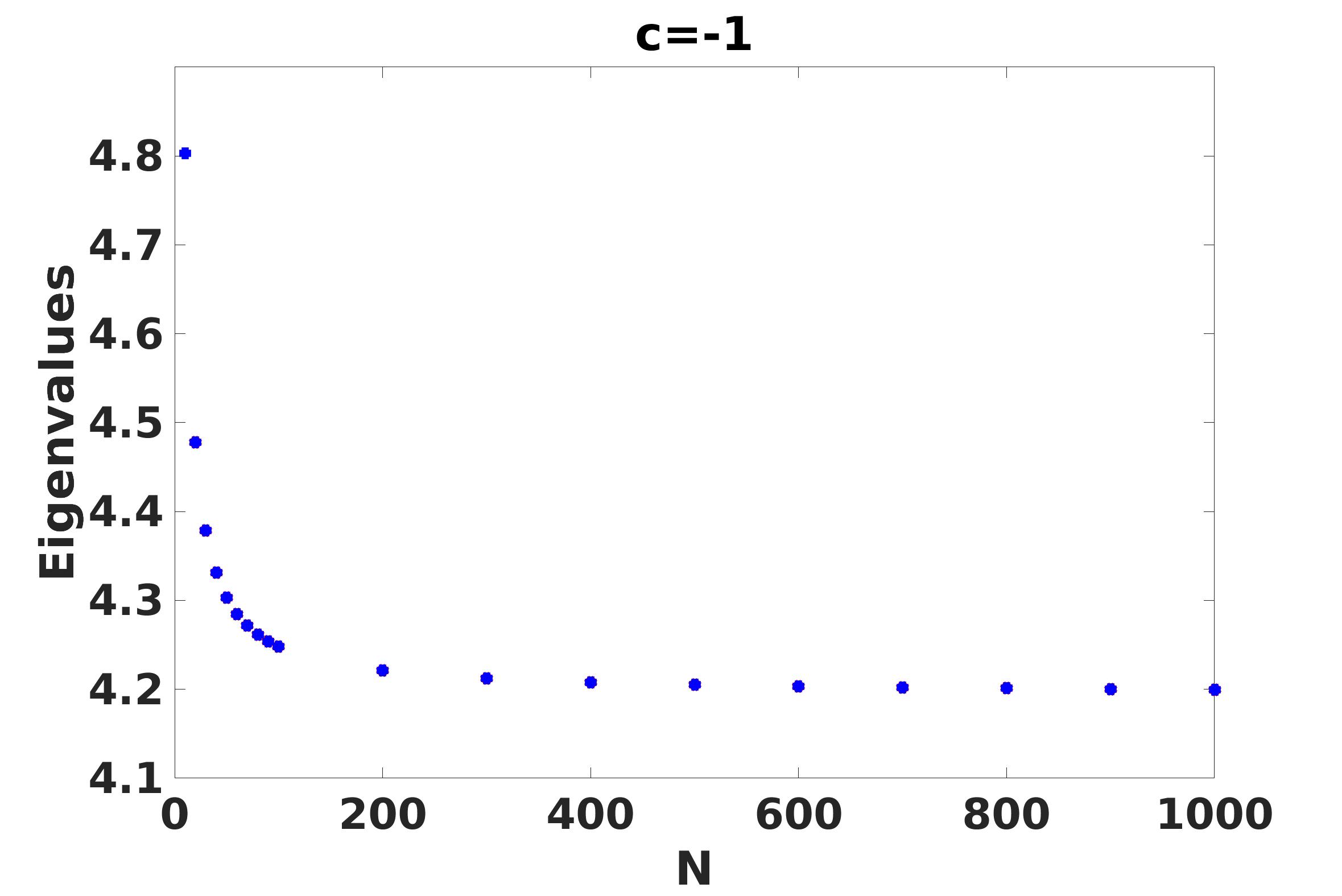}~~~~
\caption{The maximal eigenvalues of $\bar F_{\tau,h}$ for various $\alpha$ 
  and $d=2$ with $h=0.01$ as well as for $d=3$ and $h_{max}=0.2$. }
\label{fig:eigenvalues}
\end{center}
\end{figure}

The two plots in Figure \ref{fig:eigenvalues} show the largest eigenvalue of $\bar F_{\tau,h}$  for the values $\lambda=1$, $T=1$, $\alpha=1$ and $c=1$ on the right, respectively for $c=-1$ on the left, when the time step size $\tau=1/N$ varies with $N$.
They increase for $c=1$ and  decrease for $c=-1$
according to the estimate \eqref{gamma} of $\gamma$  with $N$ and converge with $N\to \infty$. Consequently they are bounded. The plots show the results for space dimension $d=2$ with spatial uniform step size $h=0.01$ (10201 grid points) as well as for $d=3$ with $h_{max}=0.2$ (1632 grid points). Visually no difference can be seen. They differ by around $1.e-12$. 
Varying the number of grid points between roughly $500$ to $40000$ leads to deviation around $1.e-12$, which reflects the independence of the spatial discretization and spatial dimension.
Also the variation of $\alpha$ has an effect of at most $1.e-9$ when considering the values 0,0.1,1,10 and 100.

As expected, the value $c$ of the zero-order term
$(cy,\eta)_{L^2(Q)} $
(or with rescaling the value of $T$) has a major effect. For $d=2, \alpha=\lambda=T=1$ and $h=0.01, N=1000$
we obtained numerically the maximal eigenvalues
listed in Table \ref{table1}.
\begin{table}[h]
\begin{center}
  \begin{tabular}{| l | r | r | r | r| r | r | r }
    \hline
    $c$             & $100$         &$10$        & $1$ &
    $-1$          & $-10$             &$-15$\\
\hline
    $\mu_{max}$&$1.00476$ &$1.04975$&$1.43204$&
   $4.24815$&$2.69623e+7$&$4.70925e+12$    \\
    \hline
  \end{tabular}
  \caption{Maximal eigenvalues of $\bar F_{\tau,h}$ for decreasing values of $c$.
  \label{table1}}
\end{center}
\end{table}
For values of $c$ lower than $-15$ the estimates are not reliable. 

The presented eigenvalue estimates indicate also bounded numbers of PCG-iterations.
As an example for the influence we consider the optimal control problem \eqref{cost1}
with
the functions
$y_0(x) = \prod_{i=1}^d \cos(x_i\pi) $ and
$y_T(x) = \prod_{i=1}^d  \cos(2x_i \pi)$
as well as 
$T=\alpha=1$,
$\lambda=10^{-3}$,
$c=-5$
unless otherwise stated.
As a starting value for the iterative solver $\bar u=0$ is employed.
One obtains a residual
$\| \lambda u_{\tau,h}+p_{\tau,h}\|_{L^2(Q)}\approx 10^{-10}$
for $d=2$
when $h=0.02$ is fixed while 
varying $N$ between 100 and 2000
with about 20-25 iterations.
When $N=50$ is fixed while varying
$h$ between $0.05$ and $0.005$ the method needs 15-17 iterations.
In three dimensions
using 4313 spatial mesh points and 
varying $N$ from 20 to 60
as well as for $N=50$ and varying the number of spatial mesh points from 508 to 11848 the iteration numbers range between  40 and 50.
Hence, in these examples the numbers of PCG-iterations stay not only bounded but also
roughly constant and quite low while varying the number of time steps
or the number of spatial discretization points.

Fixing $N=100$ and
$\dim (S_h)=  10201$ 
for $d=2$ and 
varying $c$ between 10 and $-9$ the iteration numbers are roughly between 15 and 20. 
However, for smaller values of $c$ the method stagnates and the error
$\| \lambda u_{\tau,h}+p_{\tau,h}\|_{L^2(Q)}$ increases.
For $c \leq -15$ the results are not reliable, which is due to the conditioning of the given problem.

The choice of the weight $\lambda$ for the cost $||u||^2_{L^2(Q)}$  also influences the condition number and is reflected in the iteration numbers.
For  $\lambda \in [1,10^{-1}]$  they stay constantly 6,
 then for the values $\lambda \in [10^{-1},10 ^{-5}]$ the number of iterations increases up to 132. 

%%%%%%%%%%%%%%%%%%%%%%%%
%%%%%%%%%%%%%%%%%%%%%%%%
\section{Preconditioning of the all-at-once approach}
Based on the work of Greif, Sch\"otzau  \cite{GreifSchoetzau2006} for 
preconditioning saddle point problems with highly singular (1,1) blocks we obtain the following result.
\begin{theorem}\label{GreifEtc}
Let $A_1\in \mathbb{R}^{n\times n}$ be symmetric positive semidefinite with  $\dim (ker A_1)=r$,
 $A_2\in \mathbb{R}^{m\times m}$ be a symmetric positive definite matrix and let 
$B_1 \in \mathbb{R}^{m\times n}$ have full rank, 
$B_2 \in \mathbb{R}^{m\times m}$ be invertible and $\lambda> 0$. 
Then, defining for an invertible matrix
\begin{eqnarray}
	{\cal A}= \left( 
	\begin{array}{ccc}
		A_1& 0& B_1^T\\
		0&\lambda A_2&B_2^T\\
		B_1&B_2&0
		\end{array}
	\right)
\end{eqnarray}
the preconditioner
\begin{eqnarray}
{\cal P}:= \left( 
\begin{array}{ccc}
		A_1+ B_1^T B_2^{-T} A_2 B_2^{-1} B_1& \;
		B_1^T B_2^{-T} A_2 &0\\
		A_2 B_2^{-1}B_1 &(1+\lambda) A_2 &0\\
		0&0& B_2 A_2^{-1} B_2^T
\end{array}
\right) \; ,
\end{eqnarray}
${\cal P}$ is positive definite and 
${\cal P}^{-1} {\cal A}$ has 
\begin{itemize}
	\item the eigenvalue $1$ with multiplicity $n+m$
	\item the eigenvalue $-1$ with multiplicity $r$
	\item and the remaining $m-r$ eigenvalues lie in 	$(-1,-\tfrac1{1+\lambda}] $.
\end{itemize}
\end{theorem}
\begin{proof}
Applying Theorem 2.2 in \cite{GreifSchoetzau2006} the assertion concerning the eigenvalues $1$ and $-1$ are obtained 
immediately. Moreover the remaining eigenvalues 
$\mu:=-\frac{\nu}{1+\nu} \in (-1,0)$
are given by the positive, generalized eigenvalues $\nu$
\begin{eqnarray} \nonumber
	\nu \left(
  \begin{array}{c c}
		A_1&0\\0 & \lambda A_2
	\end{array}
	\right) v
	&=&
	(B_1 \; \;  B_2)^T (B_2 A_2^{-1} B_2^T)^{-1}(B_1 \; \; B_2) v\\
	&=&
	\left(
		\begin{array}{cc}
		B_1^T B_2^{-T} A_2 B_2^{-1} B_1&\;
		B_1^T B_2^{-T} A_2  \\
		A_2 B_2^{-1}B_1 & A_2
	\end{array}
	\right) v. \label{32}
\end{eqnarray}
The second equations gives
$$
v_2=- \frac{1}{1-\nu \lambda}  B_2^{-1}B_1 v_1
$$
for $\nu \neq 1/\lambda $ and with it
$$
\nu v_1^TA_1v_1 + \frac{\nu \lambda}{1-\nu \lambda}v_1^TB_1^T B_2^{-T} A_2 B_2^{-1} B_1 v_1=0.
$$
Given $A_1$ is positive semidefinite and
$B_1^T B_2^{-T} A_2 B_2^{-1} B_1$ is positive definite this can only hold if $1-\nu \lambda<0$. This shows $-1<\mu \leq  -\frac1{1+\lambda}$.
\hfill \end{proof}

The preconditioner can be decomposed to
\begin{eqnarray*}
{\cal P}= {\left( 
\begin{array}{ccc}
	I& \tfrac1{1+\lambda} B_1^T B_2^{-T}  &0\\
	0&I  &0\\
		0&0& I
\end{array}
                     \right)
                     \left( 
\begin{array}{ccc}
		A_1+  \tfrac{\lambda}{1+\lambda} B_1^T B_2^{-T} A_2 B_2^{-1} B_1& \;
	0 &0\\
		A_2 B_2^{-1}B_1 &(1+\lambda) A_2 &0\\
		0&0& B_2 A_2^{-1} B_2^T
\end{array}
                     \right).
                     }
\end{eqnarray*}
Hence for application of $\cal P$ one needs good approximate inverses of  $A_2 $, $B_2 A_2^{-1} B_2^T$ and $A_1+  \tfrac{\lambda}{1+\lambda} B_1^T B_2^{-T} A_2 B_2^{-1} B_1$.

\medskip 
This result can be applied to the unreduced systems ${\cal T} $ and ${\cal T}_{disc} $ given in \eqref{KKTsym} and  \eqref{KKTdisc}, where in addition $\nu = 1/\lambda $ can be excluded in \eqref{32} given
$ {\cal B}_1 = {\cal K}$ is invertible.
 We obtain with  the preconditioners
\begin{align*}
  {\cal P}_{\cal T}&:=
                        {
             \left(   \begin{array}{ccc}
 		I &-\tfrac1{1+\lambda} {\cal K}^T  {\cal M}^{-1} & 0  \\
 		0& I& 0\\
 		0& 0 & I
 	\end{array}\right)
 	\left(   \begin{array}{ccc}
 e_N e_N^T \otimes M_N       
 		+ \tfrac{\lambda}{1+\lambda} {{\cal K}^T D {\cal M}^{-1}{\cal K}} &0 &0\\  
 		- D {\cal K}& (1+\lambda) D{ \cal M} & 0
 		\\
 		0& 0 & D{\cal M}
 	\end{array}\right) 
                       }
 \end{align*}
 for ${\cal T} $ and with the preconditioner
  \begin{align*}
   {\cal P}_{{\cal T}_{disc}}:=
             {
             \left(   \begin{array}{ccc}
 		I &-\tfrac1{1+\lambda} {\cal K}^T  {\cal M}^{-1} & 0  \\
 		0& I& 0\\
 		0& 0 & I
 	\end{array}\right)
 	\left(   \begin{array}{ccc}
 e_N e_N^T \otimes M_N       
 		+ \tfrac{\lambda}{1+\lambda} {{\cal K}^T D {\cal M}^{-1}{\cal K}} &0 &0\\  
 		- D {\cal K}& (1+\lambda) D{ \cal M} & 0
 		\\
 		0& 0 & D^{-1}{\cal M}
 	\end{array}\right) 
                       }
  \end{align*}
 for ${\cal T}_{disc} $ immediately the following result.
\begin{corollary}\label{PreFullSystem}
Denoting by $m:= \sum_{j=1}^Ndim (S_h^{(j)})$  the number of unknowns for each variable $u_h, y_h$ and $p_h$
the preconditioned matrices
${\cal P }^{-1}_{\cal T}{\cal T}$ and
${\cal P }^{-1}_{{\cal T}_{disc}}{\cal T}_{disc}$ have
\begin{itemize}
\item the eigenvalue 1 with multiplicity $2m $,
\item the eigenvalue -1 with multiplicity $m- dim (S_h^{(N)})$,
\item and the remaining $dim (S_h^{(N)})$ eigenvalues lie in $(  -1,-\tfrac{1}{1+\lambda}  )$.
\end{itemize}
\end{corollary}

Hence applying e.g. MINRES provides iteration numbers bounded independently of the number of grid points in time and space.

\medskip
The preconditioners are applied only approximately.
Hence approximations for the inverses of 
$ D$, $ { \cal M} $
and 
\begin{align} &
W:=e_N e_N^T \otimes M_N       
 		+ \tfrac{\lambda}{1+\lambda} {{\cal K}^T D {\cal M}^{-1}{\cal K}}
 \end{align}
are needed. 
 Similar to the control problem where  a target function is given over the whole time-horizon $[0,T]$ one may use the approximation (see \cite{PearsonStoll2013})
 \begin{gather}\label{prec}
 W
 \approx
 \tfrac{\lambda}{1+\lambda} {\cal K}^T D {\cal M}^{-1}{\cal K} =: \overline{W}. 
 \end{gather}
 Then fast approximate solvers only for
 the mass matrices $M_j$ and for $\cal K$ and ${\cal K}^T$  are required.
 Since  
\begin{align*}
  W&
 =  \mbox{
 $ {\cal K}^{T}
D 
\left[\tfrac{\lambda}{1+\lambda} I +
D^{-1} {\cal K}^{-T}
(e_N e_N^T \otimes M_N)    
{\cal K}^{-1} {\cal M}
  \right] ({\cal M}^{-1}{\cal K})^{}
  $}
       = {\cal K}^TD  \bar{F}_W {\cal M}^{-1}{\cal K}
\end{align*}
where $\bar{F}_W$ is defined as in \eqref{27} but using $\frac{\lambda}{1+ \lambda}$ instead of $\lambda$,
and employing the spectral estimate \eqref{eigenF} for $\bar{F}_W$,
one may motivate the approximation $\overline{W}$ of $W$ in \eqref{prec} also by approximating $\bar{F}_W$ with $\frac{\lambda}{1+ \lambda} I$.
 
\subsection{Numerical evidence}
To the same control problems with the same setup as for the reduced problem
 in Subsection \ref{evidencePcg} we
apply the presented preconditioner for the arising saddle point problems.
Again we employ the Cholesky decompositions of the in time constant mass and stiffness matrices.
Then only forward-backward
eliminations are applied in each time step,
which provides a solver for the preconditioner $\overline{W}$.
One may also use instead e.g.
a Chebyshev semi-iteration for preconditioning the mass matrix \cite{WathenRees2009}
and a multigrid approach for the stiffness matrix.
Alternatively also 
the all-at-once approach for
$\cal K$ and its block circulant preconditioning introduced in 
\cite{McDonaldPestanaWathen2018}
is possible.
 
Fixing in case of $d=2$ the equidistant spatial mesh to 10201 mesh points
one obtains 
$\| \lambda u_{\tau,h}+p_{\tau,h}\|_{L^2(Q)}\approx 10^{-9}$
with about 40 iterations
independently of the number of time steps $N$
when varying $N$ between  20 and 1000,
 i.e. also for 30 million unknowns.
 The $L^2$-residuals of the adjoint equation
 are also of the magnitude $10^{-9}$, while
 the $L^2$-residuals of the state equation are mostly better by a factor of $1/10$.
We obtain
similar results when fixing the number of time steps to 100 while varying the number of mesh points between 441 and 40401. 
Also considering $d=3$  one obtains iteration numbers bounded independently of the refinement  level.
With 4313 equally distributed spatial mesh points and varying the number of time steps between 20 and 1000  $L^2$-errors  of at most $10^{-8}$ are obtained with 50 to 60 iterations.
Varying the number of spatial mesh points between
$ 1632$ and $11848$ using $50$ time steps
the residuals are roughly $10^{-10}$ after
45-60 iterations.
This confirms
that the iteration numbers are bounded independently of the refinement level.

To consider the influence of the parameters 
$\lambda $ and $c$
we fix  $N=100$ and
$\dim (S_h)=  10201$  
for $d=2$ and $c=-5$ respectively $\lambda=0.001$.
The choice of $\lambda$ influences the iteration number only mildly. 
For  $\lambda \in [100,5\cdot 10^{-2}]$
they are between 9 and 11, i.e. nearly constant.
However, when $\lambda$ becomes small the preconditioner using $\overline { W} $ is close to singular and most likely due to rounding errors this results into less efficiency of the preconditioner.
The number of iterations increase with $\lambda  \in [5\cdot 10^{-2},4\cdot 10^{-3}]$  up to 21.
For values $\lambda \in [3\cdot 10^{-3}, 10^{-4}]$ the method stagnates after $26-110$ iterations.
Nevertheless the residual
$\| \lambda u_{\tau,h}+p_{\tau,h}\|_{L^2(Q)}$
and the $L^2$-residuals of the state equation as well as of the adjoint equation are still
at most $2\cdot 10^{-8}$.

As expected, the parameter $c$ is more crucial due to the error propagation of the underlying parabolic differential equation if $c$ is small.
When varying $c$ between 10 and $-8$
the iteration numbers increase from 14 to 47.
The $L^2$-residuals are less than $7 \cdot 10^{-9}$.
The results for $c \leq -9$ are not reliable or the method fails.

\section{Conclusions from the numerical results}
In accordance with the analytical results, the numerical solution of the reduced system as well as of the preconditioned full system yields
iteration numbers bounded independently of the refinement level in our numerical experiments.
They stay indeed low and roughly constant
when varying the number of unknowns.
Also, the influence of the parameter $\lambda$ and
of the function $c$  in the zero-order term of the PDE
(respectively of $T$) on the conditioning is apparent.

Given the focus on analyzing the conditioning of these systems, we employed exact solvers for
the state and adjoint equation,
i.e. of ${\cal K}$, ${\cal K}^T$ and $\cal M$
using Cholesky decompositions of the in our examples with respect to $j$ constant matrices $K_j$ and $M_j$. 
In this case, one PCG-iteration applied to the reduced system is cheaper than one with $\bar W$ preconditioned MINRES-iteration applied to the full system
by two applications of $\cal M$,
${\cal K}$ , ${\cal K}^T$ and one solve for
${\cal M}^{-1}$.
Consequently, 
in the presented examples the use of the reduced system is more efficient than the use of the full system and provides higher accuracy. This is observed also when the weight $\lambda $ for the control cost or the constant $c$ 
are varied although the analytical result Corollary
\ref{PreFullSystem}
with an exact preconditioner $W$ suggests otherwise.
This may be different if a more sophisticated approximation of the preconditioner $W$ is employed.

Moreover, in general the convergence of the PCG-method applied to the reduced system hinges on the accuracy of solving the equations with  $\cal K$ and ${\cal K}^T$.
Considering the full system these solutions play only a role for the preconditioner.
Hence, though the given  numerical results for linear parabolic control problems with a given
end-time-target function -presented for validation of the analytical results-
may suggest that it is more efficient to work with
the reduced system, this is still an open question.